\documentclass[12pt,reqno]{amsart}

\usepackage{amssymb, latexsym, amsthm}

\numberwithin{equation}{section}
\numberwithin{table}{section}

\newtheorem{thm}{Theorem}[section] %

\newtheorem{prop}[thm]{Proposition}

\newtheorem{lem}[thm]{Lemma}
\newtheorem{cor}[thm]{Corollary}
\newtheorem{defn}[thm]{Definition}

\theoremstyle{definition}
\newtheorem{rem}[thm]{Remark}
\newtheorem{exmpl}[thm]{Example}

\newcommand\kay{{\ell}}

\newcommand\Deltah{{\Delta_{2,3,7}}}

\def\Hur{\CQ_{\operatorname{Hur}}}

\newcommand{\IHur}{{I\!\Hur}}
\newcommand{\JHur}{{J\!\Hur}}

\DeclareMathOperator{\Elk}{\CQ_{Elk}}

\newcommand{\IElk}{{I\!\Elk}}

\DeclareMathOperator{\disc}{disc}
\renewcommand\th[1]{{${#1}^{\rm{th}}$}}
\def\Q{{\mathbb {Q}}}
\def\R{{\mathbb {R}}}

\def\F{{\mathbb {F}}}
\def\Z{{\mathbb {Z}}}
\def\O{{\mathcal{O}}}
\def\HC{{\mathcal H}}
\def\sub{\subseteq}
\def\s{\sigma}

\newcommand\subjectto{{\,|\ }}
\newcommand\M[1][d]{{\operatorname{M}_{#1}}}
\newcommand\opname[1]{{\operatorname{#1}}}

\newcommand\Aut{{\opname{Aut}}}
\newcommand\End{{\opname{End}}}

\newcommand\PSL[1][d]{{\operatorname{PSL}_{#1}}}
\newcommand\SL[1][d]{{\operatorname{SL}_{#1}}}
\newcommand\dimcol[2]{{[{#1}\!:\!{#2}]}}
\newcommand\Tref[1]{{Theorem~\ref{#1}}}

\newcommand\Rref[1]{{Remark~\ref{#1}}}
\newcommand\Cref[1]{{Corollary~\ref{#1}}}
\newcommand\Dref[1]{{Definition~\ref{#1}}}
\newcommand\Lref[1]{{Lemma~\ref{#1}}}

\newcommand\eq[1]{{(\ref{#1})}}

\newcommand\abs[2][]{\left|{#2}\right|_{#1}} 
\def\normali{{\lhd}} 
\newcommand\ideal[1]{{\left<{#1}\right>}}
\newcommand{\tr}[1][]{\operatorname{tr}_{#1}} %

\newcommand\tensor[1][{}]{{\otimes_{#1}}}
\def\isom{{\;\cong\;}} 
\newcommand{\set}[1]{{\{#1\}}}
\def\hra{{\,\hookrightarrow\,}}
\def\lam{{\lambda}}
\def\({\left(}
\def\){\right)}

\newcommand{\Tr}[1][]{\if!#1!{\operatorname{Tr}}\else{\operatorname{Tr}_{#1}^{\phantom{I}}}\fi} %
\long\def\forget#1\forgotten{} %

\def\ie {{\it i.e.\ }}
\def\eg {{\it e.g.\ }}
\def\cf {\hbox{\it cf.\ }}

\def\co{{\,{:}\,}}
\def\ra{{\rightarrow}}

\newcommand\ver[1]{\marginpar{\tiny Changed in Ver \VER}}

\newcommand\algint[2][]{\if!#1\relax O_{#2}
\else{O_{\!#2}^{\phantom{I}}}\fi}  \DeclareMathOperator{\syspi}{{\rm sys}\pi}

\def\CQ{\mathcal{Q}}
\def\genus{{\operatorname{\it g}}}
\def\quat{{D}} 

\newcommand\suchthat{{\,:\ \,}} %
\newcommand{\smat}[4]{{\(\!\!\begin{array}{cc}{#1}\!&\!{#2}\\[-0.1cm]{#3}\!&\!{#4}\end{array}\!\!\)}}
\newcommand{\Mat}[9]{{\begin{array}{ccc} {#1} \!&\! {#2} \!&\! {#3} \\[-0.1cm] {#4} \!&\! {#5} \!&\! {#6} \\[-0.1cm] {#7} \!&\! {#8} \!&\! {#9} \end{array}}}

\newcommand\got[1]{{\mathfrak{#1}}}

\DeclareMathOperator{\Cent}{Cent} %

\newcommand\third{{$3^{\mbox{rd}}$}}

\newif\ifXY\XYtrue

\ifXY
\usepackage{xy}
\xyoption{all} %
\fi 

\begin{document}

\title{Hurwitz quaternion order and arithmetic Riemann surfaces}

\def\BIU{Department of Mathematics, Bar Ilan University, Ramat Gan
52900, Israel}

\author[M.~Katz]{Mikhail G. Katz$^{1}$} %
\author[M.~Schaps]{Mary Schaps} %
\author[U.~Vishne]{Uzi Vishne$^2$} %
\address{\BIU} %
\email{\set{katzmik,mschaps,vishne}@math.biu.ac.il} %

\thanks{$^{1}$Supported by the Israel Science Foundation (grants no.\
84/03 and 1294/06) and the Binational Science Foundation (grant
2006393)}

\thanks{$^{2}$Supported by the EU research and training network
HPRN-CT-2002-00287, ISF Center of Excellence grant 1405/05, and
BSF grant no.~2004-083}


\subjclass[2000]{%
11R52; 
53C23,
16K20
}

\keywords{arithmetic lattice, Azumaya algebras, Fuchsian group,
Hurwitz group, Hurwitz order, hyperbolic surface, hyperbolic
reflection group, quaternion algebra, subgroup growth, systole}

\date{\today}


\begin{abstract}
We clarify the explicit structure of the Hurwitz quaternion order,
which is of fundamental importance in Riemann surface theory and
systolic geometry.
\end{abstract}

\maketitle \tableofcontents


\section{Congruence towers and the~$4/3$ bound}
\label{one}

Hurwitz surfaces are an important and famous family of Riemann
surfaces.

We clarify the explicit structure of the Hurwitz quaternion order,
which is of fundamental importance in Riemann surface theory and
systolic geometry.%
\footnote{In the literature, the term ``Hurwitz quaternion order'' has
been used both in the sense used in the present text, and in the sense
of the unique maximal order of Hamilton's rational quaternions.}
A Hurwitz surface~$X$ by definition attains the upper bound
of~$84(\genus_X - 1)$ for the order~$|{\rm Aut}(X)|$ of the
holomorphic automorphism group of~$X$, where~$\genus_X$ is its genus.
In \cite{KSV}, we proved a systolic bound
\begin{equation}
\label{11}
\syspi_1(X)\geq \tfrac{4}{3} \log(\genus_X)
\end{equation}
for Hurwitz surfaces in a principal congruence tower (see below).
Here the systole~$\syspi_1$ is the least length of a noncontractible
loop in~$X$.  The question of the existence of other congruence towers
of Riemann surfaces satisfying the bound \eqref{11}, remains open.

Marcel Berger's monograph \cite[ pp.~325--353]{Be6} contains a
detailed exposition of the state of systolic affairs up to '03.  More
recent developments are covered in \cite{SGT}.

While \eqref{11} can be thought of as a differential-geometric
application of {\em quaternion algebras}, such an application of {\em
Lie algebras\/} may be found in \cite{e7}.

We will give a detailed description of a specific quaternion algebra
order, which constitutes the arithmetic backbone of Hurwitz
surfaces. The existence of a quaternion algebra presentation for
Hurwitz surfaces is due to G.~Shimura \cite[p.~83]{Sh}.  An explicit
order was briefly described by N.~Elkies in~\cite{El0} and in
\cite{El}, with a slight discrepancy between the two descriptions, see
Remark~\ref{32b} below.  We have been unable to locate a more detailed
account of this important order in the literature.  The purpose of
this note is to provide such an account.

A Hurwitz group is by definition a (finite) group occurring as the
(holomorphic) automorphism group of a Hurwitz surface.  Such a group
is the quotient~$\Deltah/\Gamma$ of a pair of Fuchsian groups.
Here~$\Deltah$ is the~$(2,3,7)$ triangle group, while
$\Gamma\normali\Deltah$ is a normal subgroup.

Let~$\eta= 2\cos \tfrac{2\pi}{7}$ and~$K=\Q[\eta]$, a cubic extension of $\Q$.
A class of Hurwitz groups arise from ideals in the Hurwitz quaternion order
\[
\Hur \subset (\eta,\eta)_K,
\]
see \Dref{Ddef} for details.

Recall that for a division algebra~$D$ over a number field~$k$, the
discriminant~$\disc(D)$ is the product of the finite ramification
places of~$D$.  Let~$\CQ$ be an order of~$D$, and let~$\algint{K}$ be
the ring of algebraic integers in~$K$. By definition~$\CQ^1$ is the
group of elements of norm~$1$ in~$\CQ$, and a {\em principal
congruence subgroup\/} of~$\CQ^1$ is a subgroup of the form
\begin{equation}
\label{12}
\CQ^1(I) = \set{ x \in \CQ^1 \suchthat x\equiv 1\pmod{I \CQ} }
\end{equation}
where~$I \normali \algint{K}$.  Any subgroup containing such a
subgroup is called a {\em congruence subgroup}.

The Hurwitz order is described in Section~\ref{two}. In
Section~\ref{four} we verify its maximality.  The precise relationship
between the order and the~$(2,3,7)$ group is given in
Section~\ref{max}. In Section~\ref{six} we note that the Hurwitz order
is Azumaya, which implies that every ideal of the order is generated
by a central element, and every automorphism is inner. It also follows
that all quotient rings are matrix rings, and in Section~\ref{shesh}
we present some explicit examples of these quotients and their
associated congruence subgroups.

Recent developments in systolic geometry include \cite{AK, BB10,
Bal08, e7, BT, Be08, Bru, Bru2, Bru3, DKR, DR09, Elm10, EL, Gu09,
Gu11, KK, KK2, Ka4, KR2, KSh, NR, Par10, Ro, RS08, Sa08, Sa10, wen}.

\section{The Hurwitz order~$\Hur$}
\label{two}

Let~$\HC^2$ denote the hyperbolic plane.  Let~$\Aut (\HC^2)=
\PSL[2](\R)$ be its group of orientation-preserving isometries.
Consider the lattice~$\Deltah \subset \Aut(\HC^2)$, defined as the
even part of the group of reflections in the sides of the~$(2,3,7)$
hyperbolic triangle, \ie geodesic triangle with
angles~$\tfrac{\pi}{2}$,~$\tfrac{\pi}{3}$, and~$\tfrac{\pi}{7}$.  We
follow the concrete realization of~$\Deltah$ in terms of the group of
elements of norm one in an order of a quaternion algebra, given by
N.~Elkies in \cite[p.~39]{El0} and in \cite[Subsection~4.4]{El}.

Let~$K$ denote the real subfield of~$\Q[\rho]$, where~$\rho$ is a
primitive \th{7} root of unity.  Thus~$K = \Q[\eta]$, where the
element~$\eta = \rho+\rho^{-1}$ satisfies the relation
\begin{equation}
\label{21c}
\eta^3 + \eta^2 - 2\eta -1 = 0.
\end{equation}
Note the resulting identity
\begin{equation}
\label{21b}
    (2-\eta)^3 = 7 (\eta-1)^2.
\end{equation}
There are three embeddings of~$K$ into~$\R$, defined by sending~$\eta$
to any of the three real roots of~\eqref{21c}, namely
\[
2\cos\(\tfrac{2\pi}{7}\), 2\cos\(\tfrac{4\pi}{7}\),
2\cos\(\tfrac{6\pi}{7}\).
\]
We view the first embedding as the `natural' one~$K \hra \R$, and
denote the others by~$\s_1,\s_2 \co K \ra \R$. Notice that~$2
\cos(2\pi/7)$ is a positive root, while the other two are
negative.

{}From the minimal polynomial we have $\Tr[K/F](\eta) =
-1$. Multiplying \eq{21c} by the `conjugate' $\eta^3-\eta^2-2\eta+1$
gives $$(\eta^2)^3 - 5(\eta^2)^2 + 6 (\eta^2) -1 = 0,$$ so similarly
$\Tr[K/F](\eta^2) = 5$. {}The recursion relation 
\[
\Tr(\eta^{3+i}) = - \Tr(\eta^{2+i}) + 2\Tr(\eta^{1+i}) +\Tr(\eta^i)
\]
provides $\Tr(\eta^3) = -4$ and $\Tr(\eta^4) = 13$. Traces of
multiples in the integral basis $1,\eta,\eta^2$ then give the
discriminant
\[
\disc(K/\Q)=\abs{\Mat{3}{-1}{5}{-1}{5}{-4}{5}{-4}{13}} = 49 .
\]

By Minkowski's bound \cite[Subsection~30.3.3]{Hasse}, it follows
that every ideal class contains an ideal of
norm~$<\frac{3!}{3^3}\sqrt{49} < 2$, which proves that~$\algint{K} = \Z[\eta]$ is a
principal ideal domain. The only ramified prime is~$7\algint{K} =
\ideal{2-\eta}^3$, \cf \eqref{21b}
and using the fact that~$\eta-1=(\eta^2+2\eta)^{-1}$ is invertible in $\algint{K}$.
Note that the minimal polynomial~$f(t) = t^3+t^2-2t-1$ remains irreducible modulo~$2$, so
that
\begin{equation}
\label{23b} \algint{K}/\ideal{2} \isom \F_2[\bar{\eta}] = \F_8,
\end{equation}
the field with~$8$ elements.



\begin{defn}\label{Ddef}
We let~$\quat$ be the quaternion~$K$-algebra
\begin{equation}
\label{21}
(\eta,\eta)_K = K[i,j \subjectto i^2 = j^2 = \eta,\, ji = -ij].
\end{equation}
\end{defn}

As mentioned above, the root~$\eta > 0$ defines the natural imbedding
of~$K$ into~$\R$, and so~$D\otimes \R = \M[2](\R)$, and thus the
imbedding is unramified.  On the other hand, we have~$\s_1(\eta)<0$
and $\s_2(\eta)<0$, so the algebras~$D \otimes_{\sigma_1} \R$ and~$D
\otimes_{\sigma_2} \R$ are isomorphic to the standard Hamilton
quaternion algebra over~$\R$.

Moreover,~$D$ is unramified over all the finite places of~$K$
\cite[Prop.~7.1]{KSV}.
\begin{rem}
By the Albert-Brauer-Hasse-Noether theorem
\cite[Theorem~32.11]{Reiner},~$D$ is the only quaternion algebra
over~$K$ with this ramification data.
\end{rem}

Let~$\O \sub \quat$ be the order defined by
\[
\O = \algint{K}[i,j].
\]
Clearly, the defining relations of~$D$
serve as defining relations for~$\O$ as well:
\begin{equation}\label{defO}
\O \isom \Z[\eta][i,j \subjectto i^2 = j^2 = \eta,\, ji = -ij].
\end{equation}

Fix the element~$\tau = 1+\eta+\eta^2$, and define an element~$j'\in
D$ by setting
\[
j' = \frac{1}{2}(1+\eta i+\tau j).
\]

Notice that~$j'$ is an algebraic integer of~$\quat$, since the reduced
trace is~$1$, while the reduced norm is
\[
\frac{1}{4}(1-\eta\cdot \eta^2 - \eta\cdot \tau^2 + \eta^2 \cdot 0) =
-1-3\eta,
\]
so that both are in~$\algint{K}$. In particular, we have the relation
\begin{equation}
\label{jt}
{j'}^2 = j' + (1+3\eta).
\end{equation}
We define an order~$\Elk \subset \quat$ by setting
\begin{equation}
\label{22}
\Elk = \Z[\eta][i,j'].
\end{equation}
Finally, we define a new order~$\Hur \subset \quat$ by setting
\begin{equation}
\label{72}
\Hur = \Z[\eta][i,j,j'].
\end{equation}

\begin{rem}\label{32b}
There is a discrepancy between the descriptions of a maximal order
of~$\quat$ in \cite[p.~39]{El0} and in \cite[Subsection~4.4]{El}.
According to \cite[p.~39]{El0},~$\Z[\eta][i,j,j']$ is a maximal
order. Meanwhile, in \cite[Subsection~4.4]{El}, the maximal order
is claimed to be the order~$\Elk = \Z[\eta][i,j']$, described
as~$\Z[\eta]$-linear combinations of the elements~$1$, $i$, $j'$,
and~$ij'$, on the last line of \cite[p.~94]{El}.  The correct
answer is the former, \ie the order~\eqref{72}.
\end{rem}

We correct this minor error in \cite{El}, as follows.

\begin{lem}\label{2.6}
The order~$\Hur$ strictly contains~$\Elk = \Z[\eta][i,j']$.
\end{lem}
\begin{proof}
The identities \eq{jt} and
\begin{equation}
\label{23}
j'i = \eta^2 + i - ij'
\end{equation}
show that the module
\begin{equation}
\label{Q0prime}
\CQ_{\opname{Elk}}' = \Z[\eta] + \Z[\eta] i + \Z[\eta]j' +
\Z[\eta]ij',
\end{equation}
which is clearly contained in~$\Elk$, is closed under multiplication,
and thus equal to~$\Elk$.  Moreover, the set~$\set{1,i,j',ij'}$ is a
basis of~$D$ over~$K$, and a computation shows that
\[
j = \frac{-9+2\eta+3\eta^2}{7}+\frac{3-3\eta-\eta^2}{7} i+ \frac
{18-4\eta-6\eta^2} {7}j',
\]
with non-integral coefficients. Therefore~$j \not \in \Elk$.
\end{proof}

\section{Maximality of the order~$\Hur$}
\label{four}

\begin{thm}
The order~$\Hur$ is a maximal order of~$D$.
\end{thm}
\begin{proof}
By \Lref{2.6} it is enough to show that every order containing~$\Elk$
is contained in~$\Hur$. Let~$M \supseteq \Elk$ be an order, namely a
ring which is finite as a~$\algint{K}$-module, and let~$x \in
M$. Since~$\{1,i,j,ij\}$ is a~$K$-basis for the algebra~$D$, we can
write
\[
x = \tfrac{1}{2}(a+bi+cj+dij)
\]
for suitable~$a,b,c,d \in K$. Recall that every element of an
order satisfies a monic polynomial over~$\algint{K}$, so in
particular it has integral trace.  Since we have~$x,ix,jx,ijx \in
M$,
with traces~$a,\eta b,\eta c,-\eta^2d$, respectively, while the
element~$\eta = (\eta^2+\eta-2)^{-1}$ is invertible
in~$\algint{K}$, we conclude that, in fact,~$a,b,c,d \in
\algint{K}$. Now,
$$\tr(xj') = \frac{1}{4}\tr((a+bi+cj+dij)(1+\eta i+\tau j)) = \frac{1}{2}(a+\eta^2b +\eta\tau c) $$
and
$$\tr(xij') = \frac{1}{4}\tr((a+bi+cj+dij)i(1+\eta i+\tau j)) = \frac{1}{2}(\eta^2 a+\eta b-\eta^2 \tau d ).$$
Since these are integers, and since $\eta \tau \equiv \eta + 1$
and $\eta^3\tau \equiv 1$ modulo $2\algint{K}$, we have
that~$a\equiv \eta^2 b + (\eta+1)c$, and~$d \equiv \eta^3 a +
\eta^2 b \equiv \tau b + \eta c$.

It then follows that
\[
x - (\eta^2+2\eta+1)cj' - ((\eta^2+3\eta+1)c+b)ij' \in
\algint{K}[i,j],
\]
so that~$x \in \Hur$.
\end{proof}

\begin{rem}
Since~$K \Hur = D$, the center of~$\Hur$ is
$$\Cent(\Hur) = \Hur \cap \Cent(D) = \Hur \cap K = \algint{K}.$$
\end{rem}

While~$\O$ admits the presentation \eq{defO}, typical of symbol
algebras, it should be remarked that~$\Hur$ cannot have such a
presentation.
\begin{rem}\label{33}
There is no pair of anticommuting generators of~$\Hur$ over
$\Z[\eta]$.
\end{rem}
\begin{proof}
One can compute that~$\Hur/2\Hur$ is a~$2\times 2$ matrix algebra
\cite[Lemma~4.3]{KSV}, and in particular non-commutative; however
anticommuting generators will commute modulo~$2$.
\end{proof}

The prime~$2$ poses the only obstruction to the existence of an
anti-commuting pair of generators. Indeed, adjoining the fraction
$\frac{1}{2}$, we clearly have
\[
\Hur[\tfrac{1}{2}] = \O[\tfrac{1}{2}] =
\algint{K}[\tfrac{1}{2}][i,j\,|\,i^2 = j^2 = \eta,\, ji = -ij],
\]
and this is an Azumaya algebra over~$\algint{K}[\frac{1}{2}]$, see
Definition~\ref{62}. A presentation of~$\Hur$ is given in
Theorem~\ref{32}.

\section{The~$(2,3,7)$ group inside~$\Hur$}
\label{max}

The group of elements of norm~$1$ in the order~$\Hur$, modulo the
center~$\set{\pm 1}$, is isomorphic to the~$(2,3,7)$ group
\cite[p.~95]{El}. Indeed, Elkies gives the elements
\begin{eqnarray*}
g_2 & = & \frac{1}{\eta}ij, \\
g_3 & = & \frac{1}{2}(1+(\eta^2 - 2)j + (3-\eta^2)ij),  \\
g_7 & = & \frac{1}{2}((\tau-2)+(2-\eta^2)i+(\tau-3)ij),
\end{eqnarray*}
satisfying the relations~$g_2^2 = g_3^3 = g_7^7 = -1$ and~$g_2 = g_7
g_3$, which therefore project to generators of~$\Deltah \subset
\PSL[2](\R)$.

\begin{thm}
The Hurwitz order is generated, as an order, by the elements~$g_2$
and~$g_3$, so that we can write~$\Hur = \algint{K}[g_2,g_3]$.
\end{thm}

\begin{proof}
We have~$g_2,g_3,g_7 \in \Hur$ by the invertibility of~$\eta$ in
$\algint{K}$ and the equalities
\begin{eqnarray*}
g_3 & = & (3+6\eta-\eta^2) + (1+3\eta)i- (2+\eta^2) jj' -2 (ijj'-(1-\eta)ij),\\
g_7 & =
&(\tau+3\eta)+2(1+\eta)i-(\eta+\eta^2)jj'+(2-\tau)(ijj'-(1-\eta)ij)
.
\end{eqnarray*}
Conversely, we have the relations
\begin{eqnarray*}
i & = & (1+\eta)(g_3 g_2 - g_2 g_3),\\
j & = & (1+\eta)(1 + (\eta^2+\eta-1)g_2-  2 g_3),\\
j' & = & (1+\eta i)g_3 + (\eta^2 - 2)ij + j,
\end{eqnarray*}
proving the lemma.
\end{proof}

\begin{thm}
\label{32}
A basis for the order~$\Hur$ as a free module over~$\Z[\eta]$ is given
by the four elements~$1$,~$g_2$,~$g_3$, and~$g_2 g_3$.

The defining relations~$g_2^2 = -1$,~$g_3^2 = g_3 - 1$ and~$g_2g_3
+ g_3 g_2 = g_2 - (\eta^2 +\eta-1)$ provide a presentation of
$\Hur$ as an~$\algint{K}$-order.
\end{thm}

\begin{proof}
The module spanned by~$1,g_2,g_3,g_2g_3$ is closed under
multiplication by the relations given in the statement (which are
easily verified); thus~$\algint{K}[g_2,g_3] =
\opname{span}_{\algint{K}}\set{1,g_2,g_3,g_2g_3}$. The relations
suffice since the ring they define is a free module of rank~$4$,
which clearly project onto~$\Hur$.
\end{proof}

\begin{rem}
An alternative basis for the order~$\Hur$ as a free module
over~$\Z[\eta]$ is given by the four elements~$1$,~$i$,~$jj'$,
and~$\kay = ijj'-(1-\eta)ij$.
\end{rem}

\forget 
 \section{A module basis for~$\Hur$}
 \label{three}

 \begin{lem}\label{32} A basis for the order~$\Hur$ as a free module
 over~$\Z[\eta]$ is given by the four elements~$1$,~$i$,~$jj'$,
 and~$\kay = ijj'-(1-\eta)ij$.
 \end{lem}

 The element~$\kay$, as well as the formulas below, were found by
 elementary matrix operations over~$\Z[\eta]$, noting that~$\Hur
 \sub \frac{1}{2}\O$ is free, being a submodule of a free module
 over a principal ideal domain.

 \begin{proof}[Proof of \Lref{32}]
 The identities \eqref{23} and
 \begin{eqnarray}
 j j' & = & (\tau+2\eta)(i j'-2ij) - (\tau+\eta)(j'-2j), \label{jjprime} \\
 j'j & = & (3\eta+1) + j - j j' \nonumber
 \end{eqnarray}
 %
 show that~$\Hur$ is spanned, as an~$\algint{K}$-module,
 by~$1,i,j,ij,j'$ and~$ij'$.  Now the identities
 \[
 \begin{aligned}
 j & = (2\eta^2+\eta-1)[(1+3\eta) - 2 jj'] + \eta(\eta+1)[(1+3\eta)
 i -2 \kay], \\ %
 j' & = (8+17\eta+\eta^2) +(5+13\eta+3\eta^2)i -(4\eta+5\eta^2)jj'
 -(1+4\eta+3\eta^2)\kay,
 \\ %
 ij & = (1+4\eta+3\eta^2)(\eta^2 + i) -
 2(1+\eta)(\eta^2 jj'+\kay), \\ %
 ij' & =(3+8\eta+5\eta^2)( \eta^2 +i) %
 -(2+4\eta+\eta^2) (\eta^2j j' + \kay).
 \end{aligned}
 \]
 complete the proof of the lemma.
 \end{proof}

\forgotten

\forget 

It may be helpful to write down the matrix of the coefficients
expressing the elements~$j, j', ij, ij'$ in terms of the basis~$1,
i, jj', \kay$:
\[
\begin{pmatrix}
(5+10\eta-\eta^2) & -(3+7\eta+\eta^2) & -(2-2\eta-4\eta^2) &
-(2\eta+2\eta^2)
\\
(8+17\eta+8\eta^2) & (5+13\eta+3\eta^2) & -(4\eta+5\eta^2) &
-(1+4\eta+3\eta^2) &
\\
(1+5\eta+6\eta^2) & (1+4\eta+3\eta^2) & -(2+4\eta) & -(2+2\eta)
\\
(3+11\eta+10\eta^2) & (3+8\eta+5\eta^2) & -(3+7\eta+\eta^2) &
-(2+4\eta+\eta^2)
\end{pmatrix}
\]
\forgotten

\begin{rem}
Since~$\O$ is a free module of rank~$4$ over~$\algint{K}$, so
is~$\frac{1}{2}\O$, and~$\frac{1}{2}\O/\O$ is a~$4$-dimensional vector
space over~$\algint{K}/2\algint{K}$, which is the field of
order~$8$. Furthermore, one can check that~$\Hur/\O$ is a
two-dimensional subspace, namely~$\dimcol{\frac{1}{2}\O}{\Hur} =
\dimcol{\Hur}{\O} = 2^6$
where we are calculating the indices of the orders as abelian groups
.
\end{rem}

\section{Azumaya algebras}
\label{six}

We briefly describe a useful generalization of the class of
central simple algebras over fields, to algebras over commutative
rings.

\begin{defn}[\eg\ {\cite[Chapter~2]{Saltman}}]\label{62} %
Let~$R$ be a commutative ring. Let~$A$ be an~$R$-algebra which is
a faithful finitely generated projective~$R$-module. If the
natural map~$A \tensor[R] A^{\opname{op}} \ra \End_R(A)$ (action
by left and right multiplication) is an isomorphism, then~$A$ is
an {\emph{Azumaya algebra}} over~$R$.
\end{defn}

Suppose every non-zero prime ideal of~$R$ is maximal (which is the
case with every Dedekind domain, such as~$\algint{K} = \Z[\eta]$),
and let~$F$ denote the ring of fractions of~$R$. It is known that
if~$A$ is an~$R$-algebra, which is a finite module, such that
\begin{enumerate}
\item for every maximal ideal~$M \normali R$,~$A/MA$ is a central
simple algebra, of fixed degree, over~
$R/MR$; and
\item~$A\tensor[R] F$ is central simple, of the same degree, over
$F$,
\end{enumerate}
then~$A$ is Azumaya over~$R$ \cite[Theorem~2.2.a]{Saltman}. The
second condition clearly holds for~$\Hur$ over~$\algint{K}$
since~$\Hur \tensor[\algint{K}] K \isom D$.

In \cite[Lemma~4.3]{KSV} we proved the following theorem.
\newcommand\fp{{\mathfrak p}}
\begin{thm}\label{6.1}
For every ideal~$I \normali \algint{K}$, we have an isomorphism
\[
\Hur/I\!\Hur \isom \M[2] (\algint{K}/I) .
\]
\end{thm}
This was proved in \cite{KSV} for an arbitrary maximal order in a
division algebra with no finite ramification places, by
decomposing~$I$ as a product of prime power ideals, applying the
isomorphism~$\CQ/\fp^t = \CQ_{\fp} / \fp^t \CQ_{\fp}$
\cite[Section~5]{Reiner} for~$\CQ = \Hur$ ($\CQ_\fp$ being the
completion with respect to the~$\fp$-adic valuation), and using
the structure of maximal orders over a local ring
\cite[Section~17]{Reiner}.

We therefore obtain the following corollary.
\begin{cor}
The order~$\Hur$ is an Azumaya algebra.
\end{cor}

This fact has various ring-theoretic consequences. In particular,
there is a one-to-one correspondence between two-sided ideals of
$\Hur$ and ideals of its center,~$\algint{K}$
\cite[Proposition~2.5.b]{Saltman}.  Since~$\algint{K}$ is a
principal ideal domain, if follows that every two-sided ideal of
$\Hur$ is
generated by a single central element. %

Another property of Azumaya algebras is that every automorphism is
inner (namely, induced by conjugation by an invertible element)
\cite[Theorem~2.10]{Saltman}, in the spirit of the Skolem-Noether
theorem, \cf \cite[p.~107]{MR}.

\section{Quotients of~$\Hur$}\label{shesh}

The examples discussed in this section have not appeared in an
explicit form in the published literature.

In order to make \Tref{6.1} explicit, suppose~$I \normali
\algint{K}$ is an odd ideal (namely~$I+2\algint{K} = \algint{K}$).
By the inclusion~$2 \Hur \sub \O$, we have that~$\O + \IHur =
\Hur$ and~$\O \cap \IHur = I \O$.  Therefore
\[
\O/I \O = \O/(\O \cap \IHur) \isom (\IHur + \O)/\IHur = \Hur / \IHur,
\]
and so~$\O /I\O \isom \M[2](L)$ for~$L = \algint{K}/I$. {}From the
presentation of~$\O$, see \eq{defO}, it follows that
\[
\O /I\O \isom L[i,j \subjectto i^2 = j^2 = \eta, ji = -ij],
\]
which allows for an explicit isomorphism~$\O/I\O \ra \M[2](L)$.

\begin{exmpl}[First Hurwitz triplet]
\label{65}
The quotient~$\Hur/13\Hur$ can be analyzed as follows.  Since the
minimal polynomial~$\lam^3+\lam^2 - 2\lam-1$ factors over~$\F_{13}$ as
$(\lam-7)(\lam-8)(\lam-10)$, we obtain the ideal decomposition
\begin{equation}
\label{61b}
13\algint{K} = \ideal{13,\eta-7}\ideal{13,\eta-8}\ideal{13,\eta-10},
\end{equation}
and the isomorphism~$\algint{K}/\ideal{13} \ra \F_{13}\times
\F_{13}\times \F_{13}$, defined
by ~$\eta \mapsto (7,8,10)$. In fact, one has %
$$13 = \eta(\eta+2)(2\eta-1)(3-2\eta)(\eta+3),$$ where~$\eta(\eta+2)$
is invertible, and the other factors generate the ideals given
above, in the respective order; therefore, \eqref{61b} can be
rewritten as
\[
13\algint{K} = (2\eta-1) \algint{K} \cdot (3-2\eta) \algint{K} \cdot
(\eta+3) \algint{K} .
\]
An embedding~$\algint{K}[i]/ \ideal{13} \hra \M[2](\F_{13})\times
\M[2](\F_{13})\times \M[2](\F_{13})$ is obtained by mapping the
generator~$i$ via
\[
i \mapsto
\left(\smat{0}{1}{7}{0},\smat{0}{1}{8}{0},\smat{0}{1}{10}{0}
\right),
\]
satisfying the defining relation~$i^2 = \eta$. In order to extend this
embedding to~$\Hur/13\Hur$, we need to find in each case a matrix
which anti-commutes with~$i$, and whose square is~$\eta$. Namely, we
seek a matrix
\[
\smat{a}{b}{-\eta b}{-a},
\]
such that~$a^2 - \eta b^2 = \eta$ ($\eta$ stands for~$7,8$ or
$10$, respectively). Solving this equation in each case, the
map
\[
\Hur/13\Hur \ra \M[2](\F_{13})\times \M[2](\F_{13}) \times
\M[2](\F_{13})
\]
may then be defined as follows:
\[
j \mapsto \left( \smat{1}{1}{6}{12}, \smat{4}{1}{5}{9},
\smat{6}{0}{0}{7} \right).
\]

The map is obviously onto each of the components, and thus, by the
Chinese remainder theorem, onto on the product.

The three prime ideals define a triplet of principal congruence
subgroups, as in \eqref{12}. One therefore obtains a triplet of
distinct Hurwitz surfaces of genus~$14$. All three differ both in
the value of the systole and in the number of systolic loops
\cite{Vo03}.
\end{exmpl}

\begin{exmpl}[Klein quartic]
\label{6.2}
Consider the ramified prime,~$p = 7$. The minimal polynomial of~$\eta$
factors modulo~$7$ as~$(t-2)^3$, and so~$7\algint{K} = \got{p}^3$
for~$\got{p} = \ideal{\eta-2}$, see identity \eqref{21b}.  The
quotient
\[
L = \algint{K}/\got{p}^3 \isom \F_7[\epsilon|\epsilon^3 = 0]
\]
is in this case a local ring, with~$\algint{K}/7 \algint{K} \isom
L$ via~$\eta \mapsto 2+\epsilon$. The isomorphism~$\Hur/7\Hur \ra
\M[2](L)$ can be defined by~$i \mapsto \smat{0}{1}{2+\epsilon}{0}$
and~$j \mapsto (3-\epsilon+\epsilon^2) \smat{1}{0}{0}{-1}$, taking
advantage of the square root~$(3-\epsilon +\epsilon^2)^2 =
2+\epsilon$ in~$L$.

The Hurwitz surface defined by the principal congruence subgroup
associated with the ideal~$\ideal{2-\eta}$ is the famous Klein quartic, a
Hurwitz surface of genus~$3$.
\end{exmpl}

\forget 
\begin{prop} Let~$\got{p}$ be an odd prime ideal
of~$\algint{K}$, and~$I = \got{p}^t$ where~$t \geq 1$ is
arbitrary. Let~$L = \algint{K}/I$ be the quotient ring. Then
$$\Hur/\IHur \isom \O/I\O \isom \M[2](L).$$
\end{prop}
\begin{proof}
By its presentation, the ring~$\O/I\O$ is clearly isomorphic to
\[
L[i,j \subjectto i^2 = j^2 = \eta, ji = -ij],
\]
where by abuse of notation we denote the image of~$\eta$ in~$L$ by the
same letter, and it clearly remains a unit. First assume~$I$ is prime;
then~$L$ is a field, and the above algebra is well known to be a
central simple algebra of degree~$2$ over~$L$; however by Wedderburn's
theorem that there are no finite non-commutative division algebras,
such an algebra is necessarily isomorphic to~$\M[2](L)$. The same
argument, slightly generalized, covers the case where~$I$ is a prime
power, where~$L$ is a finite local (commutative) ring.

The ideals~$J = I\O$ and~$T = 2\O$ are co-prime by assumption, and
$2\Hur \sub \O$ as~$2j' = 1+\eta i + \tau j$. Therefore, the lemma
asserts that~$\Hur/\IHur \isom \O/I\O$.
\end{proof}
\forgotten

\forget
\begin{rem} An explicit isomorphism can be obtained by
sending~$i \mapsto \smat{0}{1}{\eta}{0}$ and~$j \mapsto
\smat{a}{b}{-\eta b}{-a}$, where~$a,b \in L$ solve the equation
$a^2 - \eta b^2 = \eta$. (A solution always exists, see Serre's
Course in Arithmetic).
\end{rem}
\forgotten

Having computed~$\Hur/I\Hur$ for~$I$ an odd ideal, it remains to
consider the even prime,~$2$. Recall that the order~$\Elk$ was defined
in~\eqref{22}.

\begin{thm}
Let~$I = 2^t\algint{K}$ and~$L = \algint{K}/I$. Then
$$\Hur/\IHur \isom \Elk/\IElk \isom \M[2](L).$$
\end{thm}
\begin{proof}

Since~$\tau j = 2j' - 1 - \eta i$, we have that~$\tau \Hur \sub
\Elk$; it then follows that~$\Elk + \IHur = \Hur$ and~$\Elk \cap
\IHur = \IElk$. As before,
$$\Elk/\IElk \isom \Hur / \IHur,$$ and so~$\Elk /\IElk \isom \M[2](L)$ for
$L = \algint{K}/I$.

\forget 
To complete the proof, we consider the co-prime ideals~$J = \IElk$
and~$T = \tau \Elk$ of~$\Elk$. Since~$\tau j = 2j' - 1 - \eta i$,
we have that~$\tau \Hur \sub \Elk$, and \Lref{RJS}.b asserts
that~$\Hur/\IHur \isom \Elk/\IElk$. \forgotten

Let us make this isomorphism explicit. Let~$\eta \in L$ denote the
image of~$\eta$ under the projection~$\algint{K}\ra L$. By the
relations obtained in \Lref{2.6}, we have the following presentation:
$$\Elk/\IElk \isom L[x,y\,|\,x^2 = \eta,\, y^2 = 1+3\eta+y,\,yx =
\eta^2 + x - xy]$$

Let~$b \in L$ be a solution to~$b = 2+b^2$ (such a solution exists
by Hensel's lemma, or one can iterate the defining equation).
Taking~$x' = \smat{0}{1}{\eta}{0}$ and~$y' = \smat{\eta^2}{\eta
b}{\eta^2(1-b)}{1-\eta^2}$ in~$\M[2](L)$, one verifies that~$x'$
and~$y'$ satisfy the required relations, and so~$x\mapsto x'$ and
$y\mapsto y'$ define an isomorphism.
\end{proof}

\begin{exmpl}[Fricke-Macbeath curve]
Let~$I=2\algint{K}$.  Then~$L=\algint{K}/I= \F_8$ by \eqref{23b}, and
so~$\Hur/2\Hur = \M[2](\F_8)$, as was already mentioned in
\Rref{33}. The associated quotient is
\[
\CQ_{\operatorname{Hur}}^1/\CQ_{\operatorname{Hur}}^1(I) \isom
\SL[2](\F_8)=\PSL[2](\F_8),
\]
which is the automorphism group of the corresponding Hurwitz surface
of genus~$7$, called the Fricke-Macbeath curve \cite[p.~37]{El0}.
\forget Meanwhile, the associated Fuchsian group~$\Gamma(2)$ is
the quotient of~$2\Hur$ by~$\pm Id$, and the quotient
$\Delta_{2,3,7}/\Gamma(2)$, so the quotient automorphism group is
not affected. The group \forgotten
\end{exmpl}

\begin{rem}
The quotient~$\O/2\O$ is a local commutative ring, with residue
field~$\F_8$ and a radical~$J$ whose dimension over~$\F_8$ is~$3$,
satisfying~$\dim(J^2) = 1$ and~$J^3 = 0$.
\end{rem}

\begin{proof}
Since~$ji = -ij$, the elements~$i$ and~$j$ commute modulo~$2$, so
the quotient ring is commutative. Also,~$\eta^8 \equiv \eta
\pmod{2}$. Taking the defining relations of~$\O$ modulo~$2$ we
obtain~$\O/2\O = \F_8[i,j \subjectto (i-\eta^4)^2 = (j-\eta^4)^2 =
0]$, so take~$J = \F_8 \cdot (i-\eta^4) + \F_8 \cdot (j-\eta^4) +
\F_8 \cdot (i-\eta^4)(j-\eta^4)$.
\end{proof}

\forget 
\normali \algint{K}$ are co-prime ideals, then~$\IHur \cap \JHur =
\IHur \cdot J \Hur$ (indeed write~$1 = a+b$ for~$a\in I$ and~$b \in
J$, then~$x \in \IHur \cap J \Hur$ implies~$x = ax+bx \in I \JHur + J
\IHur = \IHur \cdot \JHur$). Writing~$I = \got{p}_1^{t_1} \cdots
\got{p}_s^{t_s}$ for distinct primes~$\got{p}_i$, the Chinese
remainder theorem implies that
\[
\Hur / I \isom \prod_i \Hur/\got{p}_i^{t_i} \isom \prod_i
\M[2](\algint{K}/\got{p}_i^{t_i}) \isom \M[2](\prod_i
\algint{K}/\got{p}_i^{t_i}) \isom \M[2](I).
\]
\forgotten


\bibliographystyle{amsalpha}

\end{document}